\documentclass[a4paper]{amsart}

\usepackage[matrix,arrow,curve,tips]{xy}
            \SelectTips{cm}{}

\newtheorem{dfn}{Definition}[section]
\newtheorem{prop}[dfn]{Proposition}
\newtheorem{theo}[dfn]{Theorem}
\newtheorem{ex}[dfn]{Example}
\newtheorem{lem}[dfn]{Lemma}

\newcommand{\RR}{\mathbb{R}}
\newcommand{\CC}{\mathbb{C}}
\newcommand{\FF}{\mathbb{F}}

\newcommand{\cM}{\mathcal{M}}
\newcommand{\cN}{\mathcal{N}}

\newcommand{\ra}{\rightarrow}

\newcommand{\com}{\mathbin{{\scriptstyle \circ }}}
\newcommand{\ten}{\mathbin{\otimes}}

\newcommand{\id}{\mathord{\mathrm{id}}}
\newcommand{\pr}{\mathord{\mathrm{pr}}}
\newcommand{\supp}{\mathord{\mathrm{supp}}}
\newcommand{\uni}{\mathord{\mathrm{uni}}}
\newcommand{\inv}{\mathord{\mathrm{inv}}}
\newcommand{\mlt}{\mathord{\mathrm{mlt}}}
\newcommand{\cu}{\mathord{\epsilon}}
\newcommand{\cm}{\mathord{\Delta}}
\newcommand{\Cc}{\mathord{\mathcal{C}^{\infty}_{c}}}
\newcommand{\Cs}[1]{\mathord{\mathcal{C}_{#1}^{\infty}}}
\newcommand{\Esp}{\mathord{\mathcal{E}_{\mathit{sp}}}}
\newcommand{\Gsp}{\mathord{\mathcal{G}_{\mathit{sp}}}}
\newcommand{\pisp}{\mathord{\pi_{\mathit{sp}}}}
\newcommand{\EtGPD}{\mathsf{EtGPD}}
\newcommand{\GPD}{\mathsf{GPD}}
\newcommand{\BiALGD}{\mathsf{BiALGD}}
\newcommand{\PrBiALGD}{\mathsf{PrBiALGD}}
\newcommand{\HoALGD}{\mathsf{HoALGD}}
\newcommand{\PrHoALGD}{\mathsf{PrHoALGD}}
\newcommand{\LgHoALGD}{\mathsf{LgHoALGD}}

\title[Equivalence between the Morita categories]
      {Equivalence between the Morita categories of \'{e}tale Lie groupoids
       and of locally grouplike Hopf algebroids}

\author{J. Kali\v{s}nik}
\address{Institute of Mathematics, Physics and Mechanics,
         University of Ljubljana, Jadranska 19,
         1000 Ljubljana, Slovenia}
\email{jure.kalisnik@fmf.uni-lj.si}

\author{J. Mr\v{c}un}
\address{Department of Mathematics, University of Ljubljana,
         Jadranska 19, 1000 Ljubljana, Slovenia}
\email{janez.mrcun@fmf.uni-lj.si}

\thanks{This work was supported in part by
        the Slovenian Ministry of Science}

\begin{document}

\begin{abstract}
Any \'{e}tale Lie groupoid $G$ is completely determined by its
associated convolution algebra $\Cc(G)$ equipped with the natural
Hopf algebroid structure. We extend this result to the
generalized morphisms between \'{e}tale Lie groupoids:
we show that any
principal $H$-bundle $P$ over $G$ is
uniquely determined by the associated $\Cc(G)$-$\Cc(H)$-bimodule
$\Cc(P)$ equipped with the natural coalgebra structure.
Furthermore, we prove that the functor $\Cc$ gives
an equivalence between the Morita category
of \'{e}tale Lie groupoids and
the Morita category of locally grouplike
Hopf algebroids.
\end{abstract}

\maketitle

\section{Introduction}

The ideas and tools of noncommutative geometry have given us an
insight into a large new class of spaces, which seemed
unattainable from the point of view of the classical
topology and geometry. Lie groupoids and their convolution
algebras provide models for many such singular spaces, for example
orbifolds, spaces of orbits of Lie group actions and leaf spaces
of foliations \cite{CdSW,Con94,Hae01,Mo,MoMr}. A singular space may,
however, be represented by different Lie groupoids which are
weakly equivalent to each other. For example, the foliation
groupoids (and in particular the holonomy groupoids of foliations)
may be represented by \'{e}tale Lie groupoids
\cite{CraMoe01,MoMr}. It turns out that two Lie groupoids are
weakly equivalent if and only if they are isomorphic in the Morita
category of Lie groupoids, the category in which morphisms are
isomorphism classes of principal bundles
\cite{Hae58,HilSka,Moe91,MoMr2,Mrc99}.
We are therefore primarily interested in those algebraic
invariants of Lie groupoids which are functorially defined on the
Morita category, thus respecting the weak equivalence.

The Connes convolution algebra $\Cc(G)$ of smooth functions with
compact support \cite{Con94,Ren} on an \'{e}tale Lie groupoid $G$
is an example of
such an invariant. Indeed, the map $\Cc$ can be extended to a
functor from the Morita category of \'{e}tale Lie groupoids to the
Morita category of algebras \cite{Mrc99}. More precisely, if $G$ and
$H$ are \'{e}tale Lie groupoids and if $P$ is a principal
$H$-bundle over $G$, then the space $\Cc(P)$ of smooth functions
with compact support on $P$ has a natural structure of a
$\Cc(G)$-$\Cc(H)$-bimodule. Furthermore, the composition of
principal bundles is reflected as the tensor product of the
corresponding bimodules.

The convolution algebra $\Cc(G)$ admits an additional structure of
a coalgebra over the commutative
ring $\Cc(M)$ of smooth functions with
compact support on the base manifold $M$ of objects of $G$, which
turns $\Cc(G)$ into a Hopf algebroid over $\Cc(M)$ \cite{Mrc01,Mrc07}.
Moreover, the $\Cc(G)$-$\Cc(H)$-bimodule $\Cc(P)$ has a natural
coalgebra structure over $\Cc(M)$ as well, compatible with the
coalgebra structures on $\Cc(G)$ and $\Cc(H)$ in a natural way
\cite{Mrc01}.
The Hopf algebroid structure on $\Cc(G)$ determines the \'{e}tale
Lie groupoid $G$ uniquely \cite{Mrc07}. In fact, one can
reconstruct $G$ out of $\Cc(G)$ as the spectral \'{e}tale groupoid
associated to the Hopf algebroid $\Cc(G)$. The Hopf algebroids
isomorphic to those associated to \'{e}tale Lie groupoids can be
characterized as those which are locally grouplike
(see \cite{Mrc07} and
Definition \ref{deflocallygrouplikealgebroid} below).

In this paper we show how to reconstruct
a principal $H$-bundle $P$ over
$G$ out of the associated
$\Cc(G)$-$\Cc(H)$-bimodule $\Cc(P)$
equipped with the natural coalgebra structure.
Moreover, we show that a
$\Cc(G)$-$\Cc(H)$-bimodule $\cM$
with a $\Cc(M)$-coalgebra structure is
isomorphic to the bimodule
$\Cc(P)$ of a principal $H$-bundle $P$ over $G$ if
and only if $\cM$ is principal and locally grouplike (see Definition
\ref{deflocallygrouplikemodule}), and that the principal bundle $P$
is uniquely determined by $\cM$ up to an isomorphism.
Furthermore, we show
that locally grouplike Hopf algebroids and
locally grouplike principal bimodules form a category
$\LgHoALGD$, which
is equivalent to the Morita category $\EtGPD$
of \'{e}tale Lie groupoids
(Theorem \ref{maintheo}).
The equivalence is given by the functor
$\Cc\!:\EtGPD\ra\LgHoALGD$.

\section{Preliminaries}

For the convenience of the reader, and to fix the notations,
we begin by summarizing some basic definitions and
results that will be used in the rest of this paper.
We shall write $\FF$ for our base field, which
can be $\RR$ or $\CC$.
Throughout the paper, all manifolds and maps between them are
assumed to be smooth. This is not essential:
the results hold true
if one replaces this by any class
of differentiability $\mathcal{C}^{k}$, $k=0,1,2,\ldots\,$.
The manifolds are not assumed to be Hausdorff.

\subsection{Lie groupoids and principal bundles}
First, we recall the notion of a Lie groupoid and
the definition of the Morita category of Lie groupoids.
For detailed exposition and many examples,
we refer the reader to one of the books
\cite{Mac,MoMr,MoMr2}
and references cited there.

A Lie groupoid over a Hausdorff manifold $M$
is given by a manifold $G$ and
a structure of a category on $G$
with objects $G_{0}=M$, in which all arrows are invertible
and all the structure maps
\[
\xymatrix{G\times^{s,t}_{G_{0}}G \ar[r]^-{\mlt} &
G \ar[r]^-{\inv} & G \ar@<2pt>[r]^{s} \ar@<-2pt>[r]_{t} &
G_{0} \ar[r]^-{\uni} & G
}
\]
are smooth, with the source map $s$ a submersion.
We allow manifold $G$ to be non-Hausdorff, but we assume that
the fibers of the source map are Hausdorff.
If $g\in G$ is any arrow with source
$s(g)=x$ and target $t(g)=y$, and $g'\in G$ is another
arrow with $s(g')=y$ and
$t(g')=y'$, then the product
$g'g=\mlt(g',g)$ is an arrow from $x$ to $y'$.
The map $\uni$ assigns to each $x\in G_{0}$ the
identity arrow $1_{x}=\uni(x)$ in $G$, and we often identify
$G_{0}$ with $\uni(G_{0})$. The map $\inv$ maps each
$g\in G$ to its inverse $g^{-1}$. We write
$G(x,y)=s^{-1}(x)\cap t^{-1}(y)$.

A left action of a Lie groupoid $G$ on a manifold $P$ along a map
$\pi\!:P\ra G_{0}$ is a map $\mu\!:G\times^{s,\pi}_{G_{0}}P\ra P$,
$(g,p)\mapsto g\cdot p$, which satisfies $\pi(g\cdot p)=t(g)$,
$1_{\pi(p)}\cdot p=p$ and $g'\cdot(g\cdot p)=(g'g)\cdot p$, for all
$g',g\in G$ and $p\in P$ with $s(g')=t(g)$ and $s(g)=\pi(p)$. We
define right actions of Lie groupoids on manifolds in a
similar way.

Let $G$ and $H$ be Lie groupoids. A principal $H$-bundle over $G$
is a manifold $P$, equipped with a left action $\mu$
of $G$ along a surjective submersion $\pi\!:P\ra G_{0}$ and a right
action $\eta$ of $H$ along $\phi\!:P\ra H_{0}$, such that
$\phi(g\cdot p)=\phi(p)$, $\pi(p\cdot h)=\pi(p)$ and
$g\cdot(p\cdot h)=(g\cdot p)\cdot h$ for every $g\in G$, $p\in P$
and $h\in H$ with $s(g)=\pi(p)$ and $\phi(p)=t(h)$,
and such that
$(\pr_{1},\eta)\!:P\times^{\phi,t}_{H_{0}}H\ra P\times^{\pi,\pi}_{G_{0}}P$
is a diffeomorphism.

A map $f\!:P\ra P'$ between principal $H$-bundles $P$ and $P'$ over
$G$ is equivariant if it satisfies
$\pi'(f(p))=\pi(p)$, $\phi'(f(p))=\phi(p)$ and $f(g\cdot p\cdot
h)=g\cdot f(p)\cdot h$, for every $g\in G$, $p\in P$ and $h\in H$
with $s(g)=\pi(p)$ and
$\phi(p)=t(h)$. Any such map is automatically a
diffeomorphism. Principal $H$-bundles $P$ and $P'$ over $G$ are
isomorphic if there exists an equivariant diffeomorphism
between them.

If $P$ is a principal $H$-bundle over $G$ and $P'$ is a principal
$K$-bundle over $H$,
for another Lie groupoid $K$,
one can construct the principal $K$-bundle
$P\ten_{H}P'$ over $G$ \cite{MoMr2,Mrc99}. It is the space
of orbits
$(P\times^{\phi,\pi'}_{H_{0}}P')/H$
with respect to the natural
action of $H$ along $\phi\com\pr_{1}$
given by $(p,p')\cdot h=(p\cdot h,h^{-1}\cdot p')$.
The actions of $G$ and
$K$ on $P\ten_{H}P'$ along $\pi\ten\pi'\!:p\ten p'\mapsto\pi(p)$
respectively $\phi\ten\phi'\!: p\ten p'\mapsto\phi'(p')$ are given
by $g\cdot(p\ten p')=(g\cdot p)\ten p'$ and $(p\ten p')\cdot
k=p\ten(p'\cdot k)$, where $p\ten p'$ denotes the orbit of
$(p,p')$ in $P\ten_{H}P'$.

As an example, any Lie groupoid $G$ can be seen as a principal
$G$-bundle over $G$ with the actions given by the groupoid
multiplication along the maps $t$ respectively $s$. This
bundle behaves as the identity for
the tensor product, up to an isomorphism.

The Morita category $\GPD$ of Lie groupoids consists
of Lie groupoids
as objects and isomorphism classes of principal bundles
as morphisms between them: a principal $H$-bundle over $G$
represents a morphism from $G$ to $H$, while
the composition of morphisms is induced by the tensor product.
The morphisms in $\GPD$ are
sometimes referred to as  Hilsum-Skandalis maps or generalised
morphisms between Lie groupoids. Two Lie groupoids are Morita
equivalent if they are isomorphic in the category $\GPD$.

A Lie groupoid is \'{e}tale if all its structure maps
are local diffeomorphisms. A bisection of an \'{e}tale Lie groupoid
$G$ is an open subset $U$ of $G$ such that both $s|_{U}$ and $t|_{U}$
are injective. Any such bisection $U$
gives a diffeomorphism $\tau_{U}\!:s(U)\ra t(U)$ by
$\tau_{U}=t|_{U}\com (s|_{U})^{-1}$.

The Morita category $\EtGPD$ of \'{e}tale Lie
groupoids is the full subcategory of $\GPD$
with \'{e}tale Lie groupoids as objects.
If $G$ and $H$ are \'{e}tale Lie groupoids and $P$ a principal
$H$-bundle over $G$, then the corresponding
map $\pi\!:P\ra G_{0}$ is
automatically a local diffeomorphism.

\subsection{The bimodule associated to a principal bundle} \label{subsection2.2}
In this subsection we review the construction of the principal bimodule
assigned to a principal bundle. Our exposition closely follows
\cite{Mrc99}, where all the work is done in the Hausdorff setting.
It turns out that essentially the same ideas also work in the
non-Hausdorff case if we use the proper notion of
a smooth function with compact support.

We first recall the definition of a smooth function
with compact support on a non-Hausdorff manifold as given in
\cite{CraMoe}. Let $P$ be a manifold and let $\Cs{P}$
denote the sheaf of germs of smooth $\FF$-valued functions on
$P$. The stalk of this sheaf at a point $p\in P$ is a commutative
algebra with identity. If $P$ is Hausdorff,
we can identify
(compactly supported) smooth
functions on $P$ with (compactly supported)
continuous sections of $\Cs{P}$,
and we denote the commutative algebra of
compactly supported smooth $\FF$-valued functions on $P$ by $\Cc(P)$.
For a general $P$, we consider the space $\Gamma_{\delta}(P,\Cs{P})$ of
not-necessarily continuous sections of the sheaf $\Cs{P}$. For
any Hausdorff open subset $U\subset P$ there is a monomorphism
$\Cc(U)\ra\Gamma_{\delta}(P,\Cs{P})$, which maps $f\in\Cc(U)$ to
the extension of $f$ to $P$ by zero.
The vector space $\Cc(P)$ of smooth
functions with compact support on $P$ is by definition the image
of the map $\bigoplus_{U}\Cc(U)\ra\Gamma_{\delta}(P,\Cs{P})$,
where $U$ runs over all (or a cover of) Hausdorff open subsets of
$P$. This definition agrees with the classical
one if $P$ is Hausdorff.
We will denote the extension of $f\in\Cc(U)$ to $P$ again by
$f\in\Cc(P)$, and identify the space $\Cc(U)$ with its image in
$\Cc(P)$.

For any $f\in\Cc(P)$ we define the support of $f$ by
$\supp(f)=\{p\in P\,|\,f_{p}\neq 0\}$, where $f_{p}$ is the value
of the section $f$ at the point $p$. So defined support agrees
with the classical one in the Hausdorff case, and is a
compact set, although not closed in general. Every $f\in\Cc(P)$
with support in some Hausdorff open subset $U\subset P$ can be
identified with $f\in\Cc(U)$. We will often work with smooth
functions on the total space $P$ of a local diffeomorphism
$\pi\!:P\ra M$ into a Hausdorff manifold $M$. In this case,
we say that
$U\subset P$ is $\pi$-elementary if $\pi|_{U}$ is injective.
The open $\pi$-elementary subsets of $P$
are Hausdorff and together
they cover $P$.
We have a natural identification $\Cc(U)\cong\Cc(\pi(U))$, for
every $\pi$-elementary open subset $U\subset P$, identifying
$f_{0}\in\Cc(\pi(U))$ with
$f=f_{0}\com\pi|_{U}\in\Cc(U)\subset\Cc(P)$.

Let $\psi\!:P\ra P'$ be a smooth map. For any $p\in P$ we have a
homomorphism of algebras $\psi_{p}^{\ast}\!:(\Cs{P'})_{\psi(p)}\ra
(\Cs{P})_{p}$ given by the composition with $\psi$. If
$\psi$ is a local diffeomorphism, this homomorphism has the
inverse $\psi_{\ast p}=(\psi_{p}^{\ast})^{-1}$,
and we can define a linear map $\psi_{+}\!:\Cc(P)\ra\Cc(P')$ by
\[ \psi_{+}(f)_{p'}=
   \!\!\!\sum_{p\in\psi^{-1}(p')}\psi_{\ast p}(f_{p})\;.\]
In this way $\Cc$ becomes a functor from
the category of smooth manifolds and local diffeomorphisms
between them to the category of vector spaces.

We can use this definition of smooth functions with compact
support to construct the bimodule associated to a
principal bundle
over \'{e}tale Lie groupoids \cite{Mrc99}. Let $P$ be a
principal $H$-bundle over $G$ and let $P'$ be a principal
$K$-bundle over $H$, where $G$, $H$ and $K$ are \'{e}tale Lie
groupoids. Define a bilinear map
\[ \rho=\rho_{P,P'}\!:\Cc(P)\times\Cc(P')
   \ra\Cc(P\times^{\phi,\pi'}_{H_{0}}P')
\]
by
\[ \rho(f,f')_{(p,p')}=(\pr_{1}^{\ast})_{(p,p')}(f_{p}) \,
   (\pr_{2}^{\ast})_{(p,p')}(f'_{p'})\;,
\]
where $\pr_{1}$ and $\pr_{2}$
are the projections from $P\times_{H_{0}}P'$ to $P$ respectively $P'$.
To show that $\rho(f,f')$ is indeed a smooth function with compact
support on $P\times_{H_{0}}P'$, we can assume that $f\in\Cc(U)$
and $f'\in\Cc(U')$, where $U$ is a $\pi$-elementary open subset of
$P$ and $U'$ a $\pi'$-elementary open subset of $P'$.
The support of $\rho(f,f')$ is then contained
in the Hausdorff $\pr_{1}$-elementary open subset
$U\times_{H_{0}}U'$ of $P\times_{H_{0}}P'$, and
\begin{equation}\label{eq1}
\rho(f,f')|_{U\times_{H_{0}}U'}=(f(f'_{0}\com\phi))\com
\pr_{1}|_{U\times_{H_{0}}U'}\;,
\end{equation}
where $f'_{0}\in\Cc(\pi'(U'))$ is such that
$f'=f'_{0}\com\pi'|_{U'}$. The support
$S=\supp(f)\cap\phi^{-1}(\supp(f'_{0}))$ of the function
$f(f'_{0}\com\phi)\in\Cc(U)$ is compact and lies in the set
$\pr_{1}(U\times_{H_{0}}U')=U\cap\phi^{-1}(\pi'(U'))$. Indeed,
since $\supp(f'_{0})$ is compact in the Hausdorff manifold
$H_{0}$, it is closed, so $\phi^{-1}(\supp(f'_{0}))$ is closed as
well. The set $S$ is then a closed subspace of the compact space
$\supp(f)$ and therefore compact. This shows that $\rho(f,f')$ is
a smooth function with compact support inside $U\times_{H_{0}}U'$.

Define the map
\[
\wp=\wp_{P,P'}=q_{+}\com\rho\!:\Cc(P)\times\Cc(P')\ra
\Cc(P\ten_{H}P')\;,
\]
where $q$ is the quotient projection
$P\times_{H_{0}}P'\ra P\ten_{H}P'$, which is in fact a local
diffeomorphism. If we choose $f\in\Cc(U)$ and $f'\in\Cc(U')$ as in
the equation (\ref{eq1}), the function $\wp(f,f')$ has the support
in
$U\ten_{H}U'=q(U\times_{H_{0}}U')$ and is given by
\[
\wp(f,f')=(f(f'_{0}\com\phi))\com
((\pi|_{U})^{-1}\com(\pi\ten\pi')|_{U\ten_{H} U'})\;.
\]

Now choose another \'{e}tale Lie groupoid $L$
and a principal $L$-bundle $P''$ over $K$, and
observe that there is a natural diffeomorphism from $(P\ten_{H}
P')\ten_{K} P''$ to $P\ten_{H}(P'\ten_{K} P'')$ which maps $(p\ten p')\ten p''$
to $p\ten(p'\ten p'')$. Straight from the definition one can see
that for any $f\in\Cc(P)$, $f'\in\Cc(P')$ and $f''\in\Cc(P'')$ we
have the associativity law
\[
\wp(\wp(f,f'),f'')=\wp(f,\wp(f',f''))\;,
\]
using the natural identification $\Cc((P\ten_{H} P')\ten_{K}P'')
\cong\Cc(P\ten_{H}(P'\ten_{K} P''))$.

The equivariant diffeomorphisms
$\theta_{G,G}\!:G\ten_{G}G\ra G$, $\theta_{G,P}\!:G\ten_{G}P\ra P$ and
$\theta_{P,H}\!:P\ten_{H}H\ra P$, given by $\theta_{G,G}(g'\ten g)=g'g$,
$\theta_{G,P}(g\ten p)=g\cdot p$ and
$\theta_{P,H}(p,h)=p\cdot h$, induce the action maps
\[
\lambda_{G,G}=(\theta_{G,G})_{+}\com\wp\!:\Cc(G)\times\Cc(G)\ra\Cc(G)\;,
\]
\[
\lambda_{G,P}=(\theta_{G,P})_{+}\com\wp\!:\Cc(G)\times\Cc(P)\ra\Cc(P)
\]
and
\[
\lambda_{P,H}=(\theta_{P,H})_{+}\com\wp\!:\Cc(P)\times\Cc(H)\ra\Cc(P)\;.
\]
The map $\lambda_{G,G}$ is precisely the convolution product on
the algebra $\Cc(G)$ \cite{Con94,CraMoe,Mrc99,Mrc07}.
Furthermore, the maps
$\lambda_{G,P}$ and $\lambda_{P,H}$ turn $\Cc(P)$ into a
$\Cc(G)$-$\Cc(H)$-bimodule \cite{Mrc99}.
Explicitly, choose a
function $f\in\Cc(P)$ with support in a $\pi$-elementary open
set $U\subset P$, and functions $u\in\Cc(G)$ and
$v\in\Cc(H)$ with supports in bisections
$U_{u}\subset G$ respectively $U_{v}\subset H$.
Write $f=f_{0}\com\pi|_{U}$, $u=u_{0}\com t|_{U_{u}}$ and
$v=v_{0}\com t|_{U_{v}}$ for $f_{0}\in\Cc(\pi(U))$,
$u_{0}\in\Cc(t(U_{u}))$ and $v_{0}\in\Cc(t(U_{v}))$.
Then
\begin{equation}\label{eqleftaction}
uf=(u_{0}(f_{0}\com\tau_{U_{u}}^{-1}))\com\pi|_{\mu(U_{u}\times_{G_{0}}U)}
\end{equation}
and
\begin{equation}\label{eqrightaction}
fv=(f_{0}(v_{0}\com\phi\com(\pi|_{U})^{-1}))\com\pi|_{\eta(U\times_{H_{0}}U_{v})}\;.
\end{equation}

Now choose a principal $H$-bundle $P$ over $G$ and a principal
$K$-bundle $P'$ over $H$, and interpret $H$ as a principal $H$-bundle
over $H$. For $f\in\Cc(P)$, $f'\in\Cc(P')$ and $v\in\Cc(H)$ we
have
\[ \wp(fv,f')=\wp(\wp(f,v),f')=\wp(f,\wp(v,f'))
   =\wp(f,vf')\in\Cc(P\ten_{H}P')\;,\]
where we have identified $fv\in\Cc(P)$ with
$\wp(f,v)\in\Cc(P\ten_{H}H)$ and $vf'\in\Cc(P')$ with
$\wp(v,f')\in\Cc(H\ten_{H}P')$. The map $\wp$ thus
induces a homomorphism
\[
  \Omega=\Omega_{P,P'}\!:\Cc(P)\ten_{\Cc(H)}\Cc(P')\ra\Cc(P\ten_{H}P')
\]
of $\Cc(G)$-$\Cc(K)$-bimodules, which is in fact an isomorphism.
Indeed, this has been proven in \cite{Mrc99} in the Hausdorff
case, but literally the same proof applies in the general case as
well.

\subsection{Principal bimodules over Hopf algebroids}

Next, we review the notions of
a Hopf algebroid and of a principal bimodule over
Hopf algebroids, following \cite{Mrc01, Mrc07}.
Throughout this paper, we will assume that all our algebras
are over the field $\FF$ and that they are associative,
but not necessarily commutative.
Recall that an algebra $A$ has local
identities in a commutative subalgebra $A_{0}\subset A$
if for any $a_{1},\ldots,a_{k}\in A$ there exists
$a_{0}\in A_{0}$ such that $a_{0}a_{i}=a_{i}a_{0}=a_{i}$
for all $i=1,\ldots,k$.
A commutative
algebra has local identities if it has local identities
in itself.
A left module $\cM$ over
a commutative algebra $A_{0}$ with local identities
is locally $A_{0}$-unitary if for any
$m_{1},\ldots,m_{k}\in\cM$ there exists
$a_{0}\in A_{0}$ such that $a_{0}m_{i}=m_{i}$
for all $i=1,\ldots,k$. Analogously one defines the notions
of a right locally $B_{0}$-unitary $B_{0}$-module
and of a locally $A_{0}$-$B_{0}$-unitary $A_{0}$-$B_{0}$-bimodule,
for any commutative algebras $A_{0}$ and $B_{0}$
with local identities.
In particular, if $A$ is an algebra with local identities in
a commutative subalgebra $A_{0}\subset A$, then $A$ is
an $A_{0}$-$A_{0}$-unitary $A$-$A$-bimodule.
In this case we shall write
$A\ten^{rl}_{A_{0}}A$ for the tensor product of two copies of $A$
with respect to the right action of $A_{0}$ on the first factor and the left
action of $A_{0}$ on the second factor, while the notation
$A\ten_{A_{0}}^{ll}A$ will stand for
the tensor product taken with respect to the left action of $A_{0}$ on
both factors.

Suppose that $A_{0}$ is a commutative algebra with
local identities.  Recall that a left $A_{0}$-coalgebra
is a left $A_{0}$-unitary module $C$, together with
$A_{0}$-linear maps
$\cm\!:C\ra C\ten_{A_{0}}C$ (comultiplication)
and $\cu\!:C\ra A_{0}$ (counit) such that
$(\cu\ten\id)\com\cm=\id$, $(\id\ten\cu)\com\cm=\id$
and $(\cm\ten\id)\com\cm=(\id\com\cm)\com\cm$
(coassociativity).
We do not assume here that our coalgebras are necessarily
cocommutative, although our examples
will be such.
A homomorphism $\theta\!:C\ra C'$
of left $A_{0}$-coalgebras is a
homomorphism of left $A_{0}$-modules that respects the coalgebra
structures, i.e.\ $\cu=\cu'\com\theta$ and
$(\theta\ten\theta)\com\cm=\cm'\com\theta$.

An $A_{0}$-bialgebroid is an algebra $A$ such that
$A_{0}$ is a commutative subalgebra of $A$ in which $A$ has local
identities, together with a structure of a left $A_{0}$-coalgebra
on $A$ such that
\begin{enumerate}
\item [(i)]   $\cu|_{A_{0}}=\id$, $\cm|_{A_{0}}$ is the canonical
              embedding
              $A_{0}\subset A\ten_{A_{0}}^{ll}A$
              and the two right actions of $A_{0}$ on
              $A\ten_{A_{0}}^{ll}A$ coincide on $\cm(A)$,
\item [(ii)]  $\cu(ab)=\cu(a\cu(b))$ and
\item [(iii)] $\cm(ab)=\cm(a)\cm(b)$
\end{enumerate}
for any $a,b\in A$.
Note that by (i), the product of $A$ induces the
componentwise product
$\cm(A)\ten(A\ten_{A_{0}}^{ll}A)\ra A\ten_{A_{0}}^{ll}A$,
which is used in (iii).
The comultiplication
in an $A_{0}$-bialgebroid $A$
is also a homomorphism of right
$A_{0}$-modules with respect to any of the
right $A_{0}$-actions on
$A\ten_{A_{0}}^{ll}A$, and it induces
a homomorphism of left $A_{0}$-modules
$\overline{\cm}\!:A\ten_{A_{0}}^{rl}A\ra A\ten_{A_{0}}^{ll}A$
by $\overline{\cm}(a\ten b)=\cm(a)(a_{0}\ten b)$, where $a_{0}$
is any element of $A_{0}$ with $a_{0}b=b$.
An $A_{0}$-bialgebroid $A$ is principal if $\overline{\cm}$ is an isomorphism.
A homomorphism between $A_{0}$-bialgebroids
is a homomorphism of algebras which is also a homomorphism of
left $A_{0}$-coalgebras.

A Hopf $A_{0}$-algebroid is
an $A_{0}$-bialgebroid $A$,
together with an $\FF$-linear involution $S\!:A\ra A$ (antipode)
such that $S|_{A_{0}}=\id$,
$S(ab)=S(b)S(a)$ for any $a,b\in A$, and
$\mu_{A}\com(S\ten\id)\com\cm=\cu\com S$,
where $\mu_{A}\!:A\ten^{rl}_{A_{0}}A\ra A$ denotes the multiplication.
(This definition is slightly stronger than the definition of
a Hopf algebroid over $A_{0}$ given in \cite{Mrc07}, while the notion
of a principal Hopf algebroid
is slightly weaker than that of an \'{e}tale Hopf algebroid
given in \cite{Mrc01}.
Similar notions have been studied in
\cite{Lu,Mal,Tak,Xu} and more recently in
\cite{BloWei,BXW,BoSz,Kap}.)
A homomorphism between Hopf $A_{0}$-algebroids
is a homomorphism of $A_{0}$-bialgebroids which
intertwines the antipodes.

\begin{ex} \rm \label{ex111}
(1)
For any sheaf $\pi\!:P\ra M$ over a Hausdorff manifold $M$, the
space $\Cc(P)$ has a natural structure of a left $\Cc(M)$-coalgebra
\cite{Mrc02,Mrc0x}. The algebra $\Cc(M)$ acts on the space
$\Cc(P)$ by $(u_{0}f)_{p}=\pi^{\ast}_{p}((u_{0})_{\pi(p)})f_{p}$,
for any $u_{0}\in\Cc(M)$ and $f\in\Cc(P)$. The
comultiplication on
$\Cc(P)$ is given by
$\cm=\Omega_{\pi,\pi}^{-1}\com\mathrm{diag}_{+}$,
where
$\mathrm{diag}\!:P\ra P\times_{M}P$
is the diagonal map and
$\Omega_{\pi,\pi}\!:\Cc(P)\ten_{\Cc(M)}\Cc(P)\ra\Cc(P\times_{M}P)$
is the natural isomorphism
given by
$\Omega_{\pi,\pi}(f\ten f')_{(p,p')}
 =(\pr_{1})^{\ast}_{(p,p')}(f_{p}) (\pr_{2})^{\ast}_{(p,p')}(f_{p'})$
\cite{Mrc02,Mrc0x}. The
counit is $\cu=\pi_{+}$.
Explicitly, if $U\subset P$ is a $\pi$-elementary open subset,
$f_{0}\in\Cc(\pi(U))$
and $u_{0}\in\Cc(M)$, then
\[
u_{0}(f_{0}\com\pi|_{U})=(u_{0}f_{0})\com\pi|_{U}\;,
\]
\[
\cm(f_{0}\com\pi|_{U})=(f_{0}\com\pi|_{U})\ten(u'_{0}\com\pi|_{U})=
(u'_{0}\com\pi|_{U})\ten(f_{0}\com\pi|_{U})
\]
and
\[
\cu(f_{0}\com\pi|_{U})=f_{0}\;,
\]
where $u'_{0}\in\Cc(\pi(U))$ is any function which satisfies
$u'_{0}f_{0}=f_{0}$.

In particular, the convolution algebra $\Cc(G)$ of an
\'{e}tale Lie groupoid $G$
has a natural structure of a left
$\Cc(G_{0})$-coalgebra, induced by the target map.
The antipode
$S=\inv_{+}\!:\Cc(G)\ra\Cc(G)$ turns $\Cc(G)$
into a principal Hopf
$\Cc(G_{0})$-algebroid \cite{Mrc07}.

(2)
Any commutative algebra $A_{0}$ with local identities
is a principal Hopf $A_{0}$-algebroid in the trivial way.
\end{ex}

Suppose that
$A$ is an $A_{0}$-bialgebroid and
that $B$ is a $B_{0}$-bialgebroid.
A preprincipal $A$-$B$-bimodule
is a locally $A_{0}$-$B_{0}$-unitary $A$-$B$-bimodule $\cM$
such that
\begin{enumerate}
\item [(i)]
      the two right $B_{0}$-module structures on $\cM\ten_{A_{0}}\cM$ coincide
      on $\cm(\cM)$,
\item [(ii)]
      $\cu(mb)=\cu(m\cu(b))$ and $\cu(am)=\cu(a\cu(m))$,
\item [(iii)]
      $\cm(am)=\cm(a)\cm(m)$ and $\cm(mb)=\cm(m)\cm(b)$
\end{enumerate}
for any $a\in A$, $b\in B$ and $m\in\cM$.
The bimodule structure of
a preprincipal $A$-$B$-bimodule $\cM$ induces
the componentwise products
$\cm(A)\ten(\cM\ten_{A_{0}}\cM)\ra \cM\ten_{A_{0}}\cM$ and
$\cm(\cM)\ten (B\ten_{B_{0}}^{ll}B)\ra \cM\ten_{A_{0}}\cM$.
The existence of these two partially defined products gives the
meaning to the condition (iii) in the definition,
which also implies that $\cm$ is right $B_{0}$-linear.
The comultiplication in $\cM$
induces a homomorphism of $A_{0}$-$B$-bimodules
$\overline{\cm}\!:\cM\ten_{B_{0}} B\ra \cM\ten_{A_{0}}\cM$
by $\overline{\cm}(m\ten b)=\cm(m)(b_{0}\ten b)$,
where $b_{0}$
is any element of $B_{0}$ such that $b_{0}b=b$.
A principal $A$-$B$-bimodule is a preprincipal $A$-$B$-bimodule
$\cM$ such that $\cu$ is surjective and
$\overline{\cm}$ is an isomorphism.
A homomorphism of preprincipal
$A$-$B$-bimodules is a
homomorphism of $A$-$B$-bimodules which
is also a homomorphism of left $A_{0}$-coalgebras.

Note that a preprincipal $A$-$B$-bimodule
is in particular
a preprincipal $A_{0}$-$B$-bimodule as well as
a preprincipal $A$-$B_{0}$-bimodule.
If $\cM$ is a principal $A$-$B$-bimodule, then it
is also a principal $A_{0}$-$B$-bimodule.
Furthermore, any $A_{0}$-bialge\-broid $A$ is
also a preprincipal $A$-$A$-bimodule, which is principal
if and only if $A$ is principal as an
$A_{0}$-bialgebroid.

Let $\cM$ be a preprincipal
$A$-$B$-bimodule and let $\cN$ be a preprincipal $B$-$C$-bimodule,
for a $C_{0}$-bialgebroid $C$.
There is a natural
structure of a preprincipal $A$-$C$-bimodule
on the tensor product $\cM\ten_{B}\cN$ given by
\[
\cm(m\ten n)=\sum_{i,j}(m'_{i}\ten n'_{j})\ten(m''_{i}
\ten n''_{j})
\]
and
\[
\cu(m\ten n)=\cu(m\cu(n))\;,
\]
for any $m\ten n\in\cM\ten_{B}\cN$ with
$\cm(m)=\sum_{i}m'_{i}\ten m''_{i}$
and $\cm(n)=\sum_{j}n'_{j}\ten n''_{j}$ \cite{Mrc01}.
If $B$, $\cM$ and $\cN$ are all principal, then
$\cM\ten_{B}\cN$ is principal as well; this was proved in
\cite{Mrc01} in the  cocommutative case,
however the cocommutativity assumption
was not used in the proof.

We shall denote by $\BiALGD$ the Morita category of
bialgebroids: objects of $\BiALGD$ are pairs $(A,A_{0})$,
where $A$ is an $A_{0}$-bialgebroid, a morphism
from $(A,A_{0})$ to $(B,B_{0})$ in the category $\BiALGD$
is an isomorphism class of preprincipal
$A$-$B$-bimodules, while the composition is induced by
the tensor product.
The principal bialgebroids and the principal bimodules
form a subcategory $\PrBiALGD$ of $\BiALGD$.
Similarly, we have the Morita category
$\HoALGD$
of Hopf algebroids and isomorphism classes
of preprincipal bimodules as morphisms between them,
as well as its subcategory $\PrHoALGD$
of principal Hopf algebroids and isomorphism classes
of principal bimodules.

A smooth bialgebroid is a pair $(A,M)$, where
$M$ is a smooth Hausdorff manifold and
$A$ is a $\Cc(M)$-bialgebroid.
We have the Morita category
$\BiALGD^{\infty}$ of smooth bialgebroids,
in which a morphism
from $(A,M)$ to $(B,N)$
is an isomorphism class of preprincipal
$A$-$B$-bimodules, and the composition is induced by
the tensor product.
The principal smooth bialgebroids and the principal
bimodules
form a subcategory $\PrBiALGD^{\infty}$ of $\BiALGD^{\infty}$.
The natural functor
$\BiALGD^{\infty}\ra\BiALGD$,
which maps $(A,M)$ to $(A,\Cc(M))$,
is therefore fully-faithful, and the same is true
for its restriction
$\PrBiALGD^{\infty}\ra\PrBiALGD$.
Note that if $A$ is an $A_{0}$-bialgebroid and
if $A_{0}\cong \Cc(M)$ for a Hausdorff manifold $M$,
then $M$
is in fact determined uniquely up to
a canonical diffeomorphism:
any isomorphism $\Cc(M)\cong\Cc(N)$ is induced by a unique
diffeomorphism between $M$ and $N$.
Moreover,
the set $M$ can be in this case
identified with the set $\hat{A}_{0}$
of surjective multiplicative functionals on $A_{0}$,
so there is a natural structure of a smooth manifold
on $\hat{A}_{0}$ such that $A_{0}$ and $\Cc(\hat{A}_{0})$
are canonically isomorphic.

Analogously, we have the Morita category
$\HoALGD^{\infty}$
of smooth Hopf algebroids and isomorphism classes
of preprincipal bimodules as morphisms between them,
and also its subcategory $\PrHoALGD^{\infty}$
of principal smooth Hopf algebroids and isomorphism classes
of principal bimodules.
The natural functors
$\HoALGD^{\infty}\ra\HoALGD$ and
$\PrHoALGD^{\infty}\ra\PrHoALGD$
are fully-faithful.

\begin{ex} \rm \label{ex123}
(1)
Let $G$ and $H$ be \'{e}tale Lie groupoids and let $P$ be a
principal $H$-bundle over $G$. The
$\Cc(G)$-$\Cc(H)$-bimodule $\Cc(P)$
(Subsection \ref{subsection2.2}) carries a
natural structure of a principal $\Cc(G)$-$\Cc(H)$-bimodule
\cite{Mrc01}.
Indeed, the coalgebra structure of $P$ is
given by the sheaf $\pi\!:P\ra G_{0}$ (Example \ref{ex111} (1)),
while the principalness follows because $\overline{\cm}$
is induced by the diffeomorphism
$(\pr_{1},\eta)\!:P\times^{\phi,t}_{H_{0}}H\ra P\times^{\pi,\pi}_{G_{0}}P$.
Indeed, we have
$\Cc(P)\ten_{\Cc(G_{0})}\Cc(P)\cong \Cc(P\times_{G_{0}}P)$
(Example \ref{ex111} (1)). Furthermore, with
the methods of the proof of
\cite[Theorem 2.4]{Mrc99}, one can easily show that
we also have an isomorphism
$\Cc(P)\ten_{\Cc(H_{0})}\Cc(H)\cong \Cc(P\times_{H_{0}}H)$
induced by the map $\rho_{P,H}$
given in Subsection \ref{subsection2.2}.

The isomorphism
$\Omega$, described in Subsection \ref{subsection2.2},
respects the coalgebra structure and
is therefore an isomorphism of principal bimodules.
To sum up, we have a functor
\[
\Cc\!:\EtGPD\ra\PrHoALGD^{\infty}\;,
\]
which maps an \'{e}tale Lie groupoid $G$ to
the associated principal smooth
Hopf algebroid
$(\Cc(G),G_{0})$ and an isomorphism class
of a principal $H$-bundle $P$
over $G$ to the isomorphism class of the principal
$\Cc(G)$-$\Cc(H)$-bimodule $\Cc(P)$.

(2)
Let $x$ be a point of a Hausdorff manifold $M$. We have the
quotient epimorphism of commutative algebras
$\Cc(M)\ra (\Cs{M})_{x}$, mapping a function $f\in\Cc(M)$ to
its germ at $x$.
With respect to this epimorphism, the
space $(\Cs{M})_{x}$ has a natural structure of a principal
$(\Cs{M})_{x}$-$\Cc(M)$-bimodule.
Thus, if $A$ is a $\Cc(M)$-bialgebroid, if
$B$ is a $B_{0}$-bialgebroid and if $\cM$ a (pre)principal
$A$-$B$-bimodule, then
$(\Cs{M})_{x}\ten_{\Cc(M)}\cM$
is a (pre)principal $(\Cs{M})_{x}$-$B$-bimodule.
\end{ex}

\section{The Morita category of locally grouplike Hopf algebroids}

\subsection{Locally grouplike Hopf algebroids}

Suppose that $C=(C,\cm,\cu)$ is a left $\Cc(M)$-coalgebra,
for a Hausdorff manifold $M$.
Choose a point $x\in M$ and write $\Cc(M)_{x}=(\Cs{M})_{x}$.
There is the associated
local left $\Cc(M)_{x}$-coalgebra $(C_{x},\cm_{x},\cu_{x})$
at $x$
\cite{Mrc02,Mrc0x}, given as the quotient
$C_{x}=C/N_{x}(C)$ with respect
to the left $\Cc(M)$-submodule
\[
N_{x}(C)=\{c\in C\,|\, f_{0}c=0 \text{ for some }
f_{0}\in \Cc(M) \text{ with } (f_{0})_{x}=1\}
\]
of $C$, and with the induced coalgebra structure.
The equivalence class of an element
$c\in C$ in $C_{x}$ will be denoted by $c_{x}$.
Note that we have a natural isomorphism of left $\Cc(M)_{x}$-modules
$C_{x}\cong \Cc(M)_{x}\ten_{\Cc(M)}C$.

An element
$c\in C$ is weakly grouplike \cite{Mrc02,Mrc0x} if there exists
$c'\in C$ such that $\cm(c)=c\ten c'$.  We denote by
$G_{w}(C)$ the set of weakly grouplike elements of $C$.
We may also consider the set
$G(C_{x})=
\{\zeta\in C_{x}\,|\,\cm_{x}(\zeta)=\zeta\ten \zeta,\,\cu_{x}(\zeta)=1\}$
of grouplike elements of the coalgebra $C_{x}$.
The connection between the
weakly grouplike elements of $C$ and the
grouplike elements of
$C_{x}$ is as follows:
every $c\in G_{w}(C)$,
which is normalised at $x$ (i.e.\ $\cu(c)_{x}=1$), projects to a
grouplike element $c_{x}\in G(C_{x})$. Conversely, any
$\zeta\in G(C_{x})$ can be written as $\zeta=c_{x}$ for some
$c\in G_{w}(C)$ normalised at $x$ \cite{Mrc02,Mrc0x}.

A weakly grouplike element $a$ of
a Hopf $\Cc(M)$-algebroid $A$
is $S$-invariant if there
exists $a'\in A$ such that $\cm(a)=a\ten a'$ and
$\cm(S(a))=S(a')\ten S(a)$. Write $G^{S}_{w}(A)$ for the set
of $S$-invariant weakly grouplike elements of $A$.
The set of arrows of $A$ with target $y\in M$ is given by
$G^{S}(A_{y})=\{a_{y}\,|\, a\in G_{w}^{S}(A),\, \cu(a)_{y}=1 \}\subset G(A_{y})$.

For example,
the weakly grouplike elements of the left
$\Cc(M)$-coalgebra $\Cc(P)$
of a sheaf $\pi\!:P\ra M$ are precisely the functions on $P$
with compact support in a $\pi$-elementary open subset of
$P$ \cite{Mrc02,Mrc0x}.
Similarly, the $S$-invariant weakly grouplike elements of
the Hopf algebroid $\Cc(G)$, associated to
an \'{e}tale Lie groupoid $G$,
are the functions on $G$ with compact support in a
bisection of $G$ \cite{Mrc07}.

\begin{dfn}\rm\label{deflocallygrouplikealgebroid}
A {\em locally grouplike} Hopf algebroid
is a smooth Hopf algebroid $(A,M)$ such that
the $\Cc(M)_{y}$-module $A_{y}$ is
freely generated by the set of arrows $G^{S}(A_{y})$ with target $y$,
for every $y\in M$.
\end{dfn}

The smooth Hopf algebroid $(\Cc(G),G_{0})$ of an \'{e}tale Lie groupoid
$G$ is an example of a locally grouplike Hopf algebroid. In fact,
the converse is true as well: For any locally grouplike Hopf
algebroid $(A,M)$ there exists a spectral \'{e}tale Lie
groupoid $\Gsp(A)$ over $M$ such that
$\Cc(\Gsp(A))\cong A$
\cite{Mrc07}.
Indeed,
a smooth Hopf algebroid $(A,M)$ is
locally grouplike if and only if
the $S$-invariant
weakly grouplike elements normally generate $A$ and are normally
linearly independent. In particular,
any locally grouplike Hopf algebroid $(A,M)$
is cocommutative and principal; furthermore,
it satisfies
$(S\ten\id)\com\overline{\cm}\com(S\ten\id)\com\overline{\cm}=\id$
and
$G(A_{y})=G^{S}(A_{y})$, for any $y\in M$.

\begin{dfn}\rm\label{deflocallygrouplikemodule}
Let $(A,M)$ and $(B,N)$
be locally grouplike Hopf algebroids.
A principal $A$-$B$-bimodule $\cM$
is {\em locally grouplike} if the set of
grouplike elements $G(\cM_{x})$
freely generates the $\Cc(M)_{x}$-module
$\cM_{x}$, for every $x\in M$.
\end{dfn}

For locally grouplike Hopf algebroids
$(A,M)$ and $(B,N)$,
a principal $A$-$B$-bi\-module $\cM$ is locally grouplike if and
only if
the weakly grouplike elements of $\cM$ normally generate
$\cM$ and are normally linearly independent \cite{Mrc02,Mrc0x}. This
implies that locally grouplike principal bimodules are
cocommutative.

The principal bimodules associated to principal bundles are
locally grouplike \cite{Mrc02,Mrc0x}, and are in fact the only examples
of locally grouplike principal bimodules, up to an isomorphism. We
will show (see Theorem \ref{maintheo}) that locally grouplike Hopf
algebroids and locally grouplike principal bimodules
form a subcategory $\LgHoALGD$ of $\PrHoALGD^{\infty}$.
Moreover,
we will prove that the category $\LgHoALGD$
is equivalent to the Morita category $\EtGPD$ of \'{e}tale
Lie groupoids; an explicit equivalence is given by the functor
$\Cc\!:\EtGPD\ra\LgHoALGD$.

\subsection{The moment map}\label{subsection3.1}

Let $M$ and $N$ be Hausdorff manifolds, let
$A$ be a $\Cc(M)$-bialgebroid and let $B$ be a $\Cc(N)$-bialgebroid.
Suppose that $\cM$ is a preprincipal $A$-$B$-bimodule.
The bimodule $\cM$ is in particular a left $\Cc(M)$-coalgebra,
hence there is the associated spectral sheaf
\[
\pi=\pisp\!:\Esp(\cM)\ra M
\]
(see \cite{Mrc02,Mrc0x}). Its stalk
$\Esp(\cM)_{x}$ over a point $x$ is by definition the set
$G(\cM_{x})$, while the topology on $\Esp(\cM)$ is given by the
basis of $\pi$-elementary open
subsets $m_{W}=\{m_{x}\in G(\cM_{x})\,|\,x\in W\}$,
where $W$ is any open subset of $M$ and
$m\in G_{w}(\cM)$ any weakly grouplike element normalised on $W$
(i.e.\ $\cu(m)|_{W}=1$).

Suppose that $m\in G_{w}(\cM)$ is normalised on
an open subset $W$ of $M$. Define
a linear map
$T_{W,m}\!:\Cc(N)\ra\mathord{\mathcal{C}^{\infty}}(W)$
by
\[
T_{W,m}(v_{0})=\cu(mv_{0})|_{W}\;.
\]
This map is a homomorphism of algebras. To see this,
choose $m'\in\cM$ with $\cm(m)=m\ten m'$. Note that
$m=\cu(m)m'$. Since $\cm$ is a
homomorphism of right $\Cc(N)$-modules and the two right
$\Cc(N)$-actions on $\cm(\cM)$ coincide, we have
$\cm(mv_{0})=mv_{0}\ten m'$ for any $v_{0}\in\Cc(N)$, thus in
particular $mv_{0}\in G_{w}(\cM)$ and $mv_{0}=\cu(m'v_{0})m$.
Now the statement follows from the equalities
\begin{align*}
T_{W,m}(v_{0}v'_{0}) & =\cu(mv_{0}v'_{0})|_{W} & \\
                     & =\cu(\cu(m'v_{0})mv'_{0})|_{W} & \\
                     & =\cu(m'v_{0})|_{W} \cu(mv'_{0})|_{W}
                       &(\cu\text{ is }\Cc(M)\text{-linear}) \\
                     & =\cu(m)|_{W} \cu(m'v_{0})|_{W} \cu(mv'_{0})|_{W}
                       &(\cu(m)|_{W}=1) \\
                     & =\cu(\cu(m)m'v_{0})|_{W} \cu(mv'_{0})|_{W}
                       &(\cu\text{ is }\Cc(M)\text{-linear}) \\
                     & =\cu(mv_{0})|_{W} \cu(mv'_{0})|_{W}
                       &(m=\cu(m)m') \\
                     & = T_{W,m}(v_{0})  T_{W,m}(v'_{0})\;, &
\end{align*}
for any $v_{0},v'_{0}\in\Cc(N)$.

Now pick $x\in W$ and define a map $T_{W,m}^{x}\!:\Cc(N)\ra\FF$ by
\[
T_{W,m}^{x}(v_{0})=T_{W,m}(v_{0})(x)\;.
\]
This map is a nontrivial
(because $\cM$ is locally
unitary) multiplicative linear functional on the algebra
$\Cc(N)$, and therefore given
by the evaluation at a unique point
$z=\phi_{W,m}(x)\in N$.
Since this is true for any $x\in W$, we have the map
\[
\phi_{W,m}\!:W\ra N\;,
\]
uniquely determined by the property that
\begin{equation}\label{defphi1}
\cu(mv_{0})(x)=v_{0}(\phi_{W,m}(x))
\end{equation}
for any $v_{0}\in\Cc(N)$.

Recall that $\cM_{x}=\Cc(M)_{x}\ten_{\Cc(M)}\cM$
is a preprincipal
$\Cc(M)_{x}$-$B$-bimodule,
for any $x\in M$.
Combining the equation
(\ref{defphi1}) with the equalities $mv_{0}=\cu(m'v_{0})m$ and
$\cu(m'v_{0})|_{W}=\cu(mv_{0})|_{W}$ (the last follows from
$m=\cu(m)m'$), we get the equality
\begin{equation}\label{defphi}
(v_{0}\com\phi_{W,m})_{x}m_{x}=m_{x}v_{0}\;,
\end{equation}
which holds for any $x\in W$ and any $v_{0}\in\Cc(N)$.

Choose $x_{0}\in W$ and real functions
$\psi_{1},\ldots,\psi_{k}\in\Cc(N)$ such that
$(\psi_{1},\ldots,\psi_{k})|_{U}\!:U\ra\RR^{k}$
is a local chart on an open neighbourhood
$U$ of the point $z_{0}=\phi_{W,m}(x_{0})$ in $N$.
For each $i=1,\ldots,k$ and any $x\in W$ we have
$(\psi_{i}\com\phi_{W,m})(x)=\cu(m\psi_{i})(x)$, which shows
that $\psi_{i}\com\phi_{W,m}$ is a smooth function on $W$. From
this we conclude that $\phi_{W,m}$ is smooth on a
neighbourhood of any point $x_{0}\in W$.
Using the diffeomorphism $\pi|_{m_{W}}\!:m_{W}\ra W$ we get a smooth
map
\[
\phi_{m_{W},m}=\phi_{W,m}\com\pi|_{m_{W}}\!:m_{W}\ra N\;.
\]

Suppose that $n$ is another element of
$G_{w}(\cM)$, normalised on an open subset $V$
of $M$, such that
$m_{W}\cap n_{V}\neq\emptyset$. Choose any point
$m_{x}=n_{x}\in m_{W}\cap n_{V}$, for any $x\in V\cap W$.
By definition, this means that there exists
$u_{0}\in\Cc(M)$ such that $(u_{0})_{x}=1$ and
$u_{0}m=u_{0}n$.
We can find an open neighbourhood
$U\subset V\cap W$ of $x$ such that
$u_{0}|_{U}=1$.
We have $m_{x'}=(u_{0}m)_{x'}=(u_{0}n)_{x'}=n_{x'}$ for
all $x'\in U$ and consequently $m_{U}=n_{U}$. We will show that
the maps $\phi_{m_{W},m}$ and $\phi_{n_{V},n}$ agree on
the set $m_{U}=n_{U}$. Indeed, for any $x'\in U$ and any
$v_{0}\in\Cc(N)$ we have
\[
\cu(mv_{0})_{x'}=\cu_{x'}(m_{x'}v_{0})=
\cu_{x'}(n_{x'}v_{0})=\cu(nv_{0})_{x'}\;,
\]
which implies $v_{0}(\phi_{W,m}(x'))=v_{0}(\phi_{V,n}(x'))$. Since
$v_{0}\in\Cc(N)$ was arbitrary and since the algebra
$\Cc(N)$ separates the points of $N$, we conclude that the
maps $\phi_{W,m}$ and $\phi_{V,n}$ agree on $U$. The same is then
true for the maps $\phi_{m_{W},m}$ and $\phi_{n_{V},n}$ on
the set $m_{U}=n_{U}$. Gluing together these locally defined maps
we get a globally defined smooth map
\[
\phi\!:\Esp(\cM)\ra N\;,
\]
the {\em moment map} of the preprincipal $A$-$B$-bimodule $\cM$.

For any $p\in G(\cM_{x})$ we define the {\em effect germ}
$\phi^{\pi}_{p}=\phi_{p}\com (\pi_{p})^{-1}$ of $p$,
which is the germ at $x$ of
a map $(M,x)\ra (N,\phi(p))$. The equation (\ref{defphi})
yields
\[
((v_{0})_{\phi(p)}\com\phi^{\pi}_{p}) p=pv_{0}
\]
for any $v_{0}\in\Cc(N)$. In particular, this equation shows that
the left $\Cc(M)_{x}$-submodule
\[
\cM^{p}_{x}=\Cc(M)_{x}p\subset \cM_{x}
\]
is also a right $\Cc(N)$-submodule of $\cM_{x}$. Moreover, the right
$\Cc(N)$-action on $\cM^{p}_{x}$ induces a right $\Cc(N)_{\phi(p)}$-action, so
$\cM^{p}_{x}$ is in fact a $\Cc(M)_{x}$-$\Cc(N)_{\phi(p)}$-bimodule.

\subsection{Tensor product of locally grouplike principal bimodules}

Let $(A,M)$, $(B,N)$ and $(C,N')$
be locally grouplike Hopf algebroids, let
$\cM$ be a locally grouplike
principal $A$-$B$-bimodule and let
$\cN$ be a locally grouplike
principal $B$-$C$-bimodule.
Write $A_{0}=\Cc(M)$ and $B_{0}=\Cc(N)$.
We know that $\cM\ten_{B}\cN$ is a principal
$A$-$C$-bimodule. Furthermore,
since $\cM$ is also a preprincipal
$A$-$B_{0}$-bimodule and $\cN$ is a principal
$B_{0}$-$C$-bimodule, it follows that
$\cM\ten_{B_{0}}\cN$ is a preprincipal $A$-$C$-bimodule.

Take any $x\in M$. Since $(A_{0})_{x}$
is a principal Hopf $(A_{0})_{x}$-algebroid
and also a principal $(A_{0})_{x}$-$A_{0}$-bimodule
(Example \ref{ex123} (2)),
we have an isomorphism of preprincipal
$(A_{0})_{x}$-$C$-bimodules
\[
(\cM\ten_{B_{0}}\cN)_{x} \cong (A_{0})_{x}\ten_{A_{0}}\cM\ten_{B_{0}}\cN
   \cong \cM_{x}\ten_{B_{0}}\cN
\]
and an isomorphism of principal $(A_{0})_{x}$-$C$-bimodules
\[ (\cM\ten_{B}\cN)_{x} \cong (A_{0})_{x}\ten_{A_{0}}\cM\ten_{B}\cN
   \cong \cM_{x}\ten_{B}\cN\;.
\]
Furthermore,
the local coalgebra $(\cM\ten_{A_{0}}\cM)_{x}$
of the product coalgebra $\cM\ten_{A_{0}}\cM$
is simply the product coalgebra
$\cM_{x}\ten_{(A_{0})_{x}}\cM_{x}=\cM_{x}\ten_{A_{0}}\cM_{x}$
over $(A_{0})_{x}$ (see also
\cite[Proposition 2.1]{Mrc02} and \cite{Mrc0x}).

\begin{lem}\label{lemmadeltabar}
The map
$\overline{\cm}\!:\cM\ten_{B_{0}}B\ra\cM\ten_{A_{0}}\cM$
is an isomorphism of left $A_{0}$-coalgebras, and
induces an isomorphism
$\overline{\cm}_{x}\!:\cM_{x}\ten_{B_{0}}B\ra\cM_{x}\ten_{A_{0}} \cM_{x}$
of left $(A_{0})_{x}$-coalgebras, for every $x\in M$.
\end{lem}

\begin{proof}
Since the left $A_{0}$-module
$\cM$ is generated by $G_{w}(\cM)$,
a straightforward
calculation shows that $\overline{\cm}$ is an isomorphism of
coalgebras over $A_{0}$. Thus we have the induced
isomorphism of $(A_{0})_{x}$-coalgebras
$(\cM\ten_{B_{0}}B)_{x}\ra(\cM\ten_{A_{0}}\cM)_{x}$,
which we combine with the isomorphisms
$(\cM\ten_{A_{0}}\cM)_{x}\cong \cM_{x}\ten_{A_{0}}\cM_{x}$
and
$(\cM\ten_{B_{0}}B)_{x} \cong
(A_{0})_{x}\ten_{A_{0}}\cM\ten_{B_{0}}B \cong
\cM_{x}\ten_{B_{0}}\cN$
to obtain $\overline{\cm}_{x}$.
Alternatively, one can also describe
$\overline{\cm}_{x}$ as the isomorphism
$\overline{(\cm_{x})}$ associated to the
principal $(A_{0})_{x}$-$B$-bimodule
$\cM_{x}\cong (A_{0})_{x}\ten_{A_{0}}\cM$.
\end{proof}

We will next describe the grouplike elements of the
coalgebra $\cM_{x}\ten_{B_{0}}\cN$.
Since $\cM$ is a locally grouplike principal $A$-$B$-bimodule, the
set $G(\cM_{x})$ freely generates
the left $(A_{0})_{x}$-module $\cM_{x}$, for every $x\in M$.
Choose any $p\in G(\cM_{x})$.
Observe that we have
the natural isomorphisms of left $(A_{0})_{x}$-modules
\[
\cM^{p}_{x}\ten_{B_{0}}\cN \cong
\cM^{p}_{x}\ten_{(B_{0})_{\phi(p)}} (B_{0})_{\phi(p)} \ten_{B_{0}}\cN
\cong
\cM^{p}_{x}\ten_{(B_{0})_{\phi(p)}}\cN_{\phi(p)}\;.
\]
Since $\cN$ is locally grouplike as well, it follows that
$\cN_{\phi(p)}$ is a free left
$(B_{0})_{\phi(p)}$-module with the basis
$G(\cN_{\phi(p)})$. For any $q\in G(\cN_{\phi(p)})$
we have the left
$(A_{0})_{x}$-module
$\cM^{p}_{x}\ten_{(B_{0})_{\phi(p)}}\cN_{\phi(p)}^{q}
\cong (A_{0})_{x}(p\ten q)$.
Therefore the left $(A_{0})_{x}$-modules
\begin{align*}
\cM^{p}_{x}\ten_{B_{0}}\cN
&\cong
\cM^{p}_{x}\ten_{(B_{0})_{\phi(p)}}\cN_{\phi(p)}
\cong
\cM^{p}_{x}\ten_{(B_{0})_{\phi(p)}}
\!\!\!\bigoplus_{q\in G(\cN_{\phi(p)})}\!\!\!\!\!\!
\cN_{\phi(p)}^{q} \\
&\cong
\!\!\!\bigoplus_{q\in G(\cN_{\phi(p)})}\!\!\!\!\!\!
(\cM^{p}_{x}\ten_{(B_{0})_{\phi(p)}}\cN_{\phi(p)}^{q})
\end{align*}
are all free. The explicit base of the
free $(A_{0})_{x}$-module $\cM^{p}_{x}\ten_{B_{0}}\cN$
is the set $L_{p}=\{[p,q]\,|\, q\in G(\cN_{\phi(p)})\}$, where
$[p,q]$ denotes the element of
$\cM^{p}_{x}\ten_{B_{0}}\cN$ which corresponds to the element
$p\ten q\in\cM^{p}_{x}\ten_{(B_{0})_{\phi(p)}}\cN_{\phi(p)}^{q}$
by the above isomorphism.

\begin{lem}\label{lemmagrouplikeelements}
(i)
The coalgebra $\cM_{x}\ten_{B_{0}}\cN$ over $(A_{0})_{x}$ is
freely generated by its grouplike elements
$G(\cM_{x}\ten_{B_{0}}\cN)
=\{[p,q]\,|\,
p\in G(\cM_{x}),\, q\in G(\cN_{\phi(p)})\}$.

(ii)
The coalgebra $\cM_{x}\ten_{B_{0}}\cN$ is
a direct sum
$\bigoplus_{p\in G(\cM_{x})}\cM^{p}_{x}\ten_{B_{0}}\cN$
of subcoalgebras over $(A_{0})_{x}$, and each subcoalgebra
$\cM^{p}_{x}\ten_{B_{0}}\cN$ is freely generated by its grouplike
elements $G(\cM^{p}_{x}\ten_{B_{0}}\cN)
=\{[p,q]\,|\,
q\in G(\cN_{\phi(p)})\}$.

(iii) The coalgebra
$\cM_{x}\ten_{A_{0}}\cM_{x}$ is a free left $(A_{0})_{x}$-module,
generated by its grouplike elements
$G(\cM_{x}\ten_{A_{0}}\cM_{x})=
\{p\ten p'\,|\,p,p'\in G(\cM_{x})\}
\cong G(\cM_{x})\times G(\cM_{x})$.
\end{lem}

\begin{proof}
(i)
It is clear that $\cM_{x}\ten_{B_{0}}\cN$ is a free
$(A_{0})_{x}$-module with the basis $L=\cup_{p} L_{p}$.
We need to show that
$G(\cM_{x}\ten_{B_{0}}\cN)=L$.
Straight from the definition of the structure maps of the
coalgebra $\cM_{x}\ten_{B_{0}}\cN$ it follows
$L\subset G(\cM_{x}\ten_{B_{0}}\cN)$.
To prove the converse inclusion, choose any
$u\in G(\cM_{x}\ten_{B_{0}}\cN)$ and write it in the form
$u=\sum_{p, q} a_{p q}\, [p,q]$
(where $p\in G(\cM_{x})$ and $q\in G(\cN_{\phi(p)})$)
for uniquely determined $a_{p q}\in (A_{0})_{x}$.
Then we have
\[
u\ten u= \!\!\sum_{p, q,p', q'}\!\!
a_{p q}a_{p' q'}[p,q]\ten[p',q']
\]
and
\[
\cm(u)=\sum_{p, q} a_{p q}
[p,q]\ten[p,q]\;.
\]
Since the element $u$ is grouplike, we have $\cm(u)=u\ten u$ and
$\cu(u)=1$. Therefore, by checking the components of $\cm(u)$ and
$u\ten u$, we see that
\[
a_{p q}^{2}=a_{p q}\;\;\;\;\;\;\;\; \text{for all }p,\, q\;,
\]
\begin{equation}\label{eq7}
a_{p q}a_{p' q'}=0 \;\;\;\;\;\;\;\;
\text{if }p\neq p'\text{ or } q\neq q'\;,
\end{equation}
and
\begin{equation}\label{eq6}
\sum_{p, q}a_{p q}=1\;.
\end{equation}
The equation (\ref{eq6}) implies that there exist $p_{0}$ and $q_{0}$
such that $a_{p_{0} q_{0}}$ is invertible in
$(A_{0})_{x}$. Combining this with the equation
(\ref{eq7}) we see that $a_{p q}=0$ if $p\neq p_{0}$
or $q\neq q_{0}$,
and that $a_{p_{0}q_{0}}=1$.
We conclude that
$u=[p_{0},q_{0}]\in L$. Parts (ii) and (iii)
follow analogously.
\end{proof}

\begin{prop}
The locally grouplike Hopf algebroids and
the locally grouplike principal
bimodules form a subcategory $\LgHoALGD$
of the category
$\PrHoALGD^{\infty}$.
\end{prop}

\begin{proof}
Let $(A,M)$, $(B,N)$ and $(C,N')$
be locally grouplike Hopf algebroids, let
$\cM$ be a locally grouplike principal $A$-$B$-bimodule
and let $\cN$ be a locally grouplike principal $B$-$C$-bimodule.
We have to show that the tensor product $\cM\ten_{B}\cN$
is locally grouplike.
If we restrict the isomorphism
$\overline{\cm}_{x}$ (Lemma \ref{lemmadeltabar})
to the submodule $\cM^{p}_{x}\ten_{B_{0}}B$,
for some $x\in M$ and $p\in G(\cM_{x})$,
we obtain an isomorphism
$\cM^{p}_{x}\ten_{B_{0}}B\cong
\cM^{p}_{x}\ten_{A_{0}}\cM_{x}\cong \cM_{x}$
of coalgebras over
$(A_{0})_{x}$. Indeed, this follows
from Lemma \ref{lemmagrouplikeelements}
by considering the explicit basis of both
$(A_{0})_{x}$-modules. It follows that
we have the isomorphisms of coalgebras over $(A_{0})_{x}$
\[
(\cM\ten_{B}\cN)_{x}\cong
\cM_{x}\ten_{B}\cN\cong(\cM^{p}_{x}\ten_{B_{0}}B)\ten_{B}\cN\cong
\cM^{p}_{x}\ten_{B_{0}}\cN\;.
\]
Since $\cM^{p}_{x}\ten_{B_{0}}\cN$ is
freely generated by its grouplike elements by Lemma
\ref{lemmagrouplikeelements}, so is $(\cM\ten_{B}\cN)_{x}$.
\end{proof}

The category $\LgHoALGD$ will be referred to as the
{\em Morita category of locally grouplike Hopf algebroids}.

\section{The principal bundle associated to a locally grouplike principal
bimodule} \label{section3}

Let $(A,M)$ and $(B,N)$ be locally grouplike Hopf
algebroids and let $\cM$ be a locally grouplike principal
$A$-$B$-bimodule.
Denote by $G=\Gsp(A)$ and $H=\Gsp(B)$ the spectral
\'{e}tale Lie groupoids associated to $(A,M)$ respectively
$(B,N)$ \cite{Mrc07}.
In particular, we have $G_{0}=M$ and $H_{0}=N$.
Recall that there are natural
isomorphisms of Hopf algebroids $\Cc(G)\cong A$ and
$\Cc(H)\cong B$, so we may regard $\cM$ as a principal
$\Cc(G)$-$\Cc(H)$-bimodule. In this section
we will construct an associated principal
$H$-bundle $P=\cM_{\ast}$ over $G$ such that the principal
$\Cc(G)$-$\Cc(H)$-bimodules $\Cc(P)$ and $\cM$ are isomorphic.

For the manifold $P$ we take
the total space of the spectral sheaf
$\pisp\!:\Esp(\cM)\ra G_{0}$ associated to the
$\Cc(G_{0})$-coalgebra $\cM$
(see Subsection \ref{subsection3.1}),
\[
P = \Esp(\cM)\;.
\]
Put $\pi=\pisp$.
We also have the associated moment map
$\phi\!:P\ra H_{0}$ constructed in Subsection \ref{subsection3.1}.

The locally grouplike Hopf algebroid $(A,G_{0})$ is in
particular a coalgebra over $\Cc(G_{0})$. Thus we have the
associated spectral sheaf $\pisp(A)\!:\Esp(A)\ra G_{0}$,
which is in fact equal to the target map of the
spectral \'{e}tale Lie groupoid
$t\!:\Gsp(A)\ra G_{0}$
because $(A,G_{0})$ is locally grouplike.
In particular, we have $\Esp(A)=\Gsp(A)=G$, while
$G^{S}(A_{y})=G(A_{y})=t^{-1}(y)$
for any $y\in G_{0}$. Recall that
an arrow $g\in G(x,y)$ can be represented
as $g=a_{y}$
by an element
$a\in G_{w}^{S}(A)$ normalised on an open neighbourhood
$W_{a}\subset G_{0}$ of $y$. Such a pair $(W_{a},a)$
induces a
diffeomorphism $\tau_{W_{a},a}\!:V_{W_{a},a}\ra W_{a}$,
defined on an open
subset $V_{W_{a},a}\subset G_{0}$, which is determined by the
property that
\begin{equation}\label{algebroidrelation}
u_{0}a=a(u_{0}\com\tau_{W_{a},a})
\end{equation}
for any $u_{0}\in\Cc(W_{a})$ \cite{Mrc07}.
We have $x=s(a_{y})=\tau_{W_{a},a}^{-1}(y)$.
The {\em effect germ} $\tau_{g}=(\tau_{W_{a},a})_{x}$ of $g$,
which is a germ at $x$
of a diffeomorphism $(G_{0},x)\ra (G_{0},y)$,
depends only on $g$ and not
on the choice of $a$ and $W_{a}$. The equation
(\ref{algebroidrelation})
may be rewritten as
\[
((u'_{0})_{x}\com \tau_{g}^{-1})g=g u'_{0}\;,
\]
for any $u'_{0}\in\Cc(G_{0})$.
For another arrow $g'\in G(y,y')$, represented
as $g'=a'_{y'}$
by $a'\in G_{w}^{S}(A)$ normalised on an open neighbourhood
$W_{a'}\subset G_{0}$ of $y'$, the product of $g'$ and $g$
is given by
$g'g=a'_{y'}a_{y}=(a'a)_{y'}$.
Any function $a_{0}\in\Cc(G_{0})$ with
$(a_{0})_{x}=1$ represents the identity arrow $1_{x}=(a_{0})_{x}$ at $x$,
while $a_{y}^{-1}=(S(a))_{x}$.

\subsection{Construction of the actions of $G$ and $H$ on $P$}

Let $g\in G(x,y)$, $p\in P$ and $h\in H(z',z)$
be such that $\pi(p)=x$ and $\phi(p)=z$.
Choose $a\in G_w^S(A)$, $m\in G_w(\cM)$
and $b\in G_{w}^{S}(B)$ such that
$a$ is normalised on an open neighbourhood $W_{a}\subset G_{0}$ of $y$,
$m$ is normalised on an open neighbourhood $W_{m}\subset G_{0}$ of $x$
and
$b$ is normalised on an open neighbourhood $W_{b}\subset H_{0}$ of $z$
with $g=a_{y}$, $p=m_{x}$ and $h=b_{z}$.

\begin{lem}\label{lem1}
The elements $am$ and $mb$ of $\cM$ are
weakly grouplike and normalised on
$W_{am}=\tau_{W_{a},a}(V_{W_{a},a}\cap W_{m})$ respectively
$W_{mb}=\pi(m_{W_{m}}\cap\phi^{-1}(W_{b}))$.
\end{lem}

\begin{proof}
From $\cm(a)=a\ten a'$, $\cm(b)=b\ten b'$ and
$\cm(m)=m\ten m'$ we get
$\cm(am)=\cm(a)\cm(m)=(a\ten a')(m\ten m')=am\ten a'm'$
and $\cm(mb)=(m\ten m')(b\ten b')=mb\ten m'b'$.
Next, for any $y'\in W_{am}\subset W_{a}$ and any function
$u_{0}\in\Cc(W_{am})$ with $(u_{0})_{y'}=1$ we have
\begin{align*}
\cu(am)(y') & = u_{0}(y')\cu(am)(y')
              = \cu(u_{0}am)(y') & \\
            & = \cu(a(u_{0}\com\tau_{W_{a},a})m)(y')
            & (\text{by (\ref{algebroidrelation})})\\
            & = \cu(a(u_{0}\com\tau_{W_{a},a})\cu(m))(y')
            & \\
            & = \cu((u_{0}(\cu(m)\com\tau_{W_{a},a}^{-1}))a)(y')
            & (\text{by (\ref{algebroidrelation})})\\
            & =
            u_{0}(y')(\cu(m)\com\tau_{W_{a},a}^{-1})(y')\cu(a)(y')
            & \\
            & = \cu(m)(\tau_{W_{a},a}^{-1}(y'))\cu(a)(y')=1\;,
\end{align*}
which proves that $am$ is normalised on $W_{am}$.
Finally, choose any $x'\in W_{mb}\subset W_{m}$. Then $\phi(m_{x'})\in
W_{b}$ and therefore $\cu(b)(\phi(m_{x'}))=1$. The equation (\ref{defphi1})
yields
\[
\cu(mb)(x')=\cu(m\cu(b))(x')=\cu(b)(\phi(m_{x'}))=1\;,
\]
so $mb$ is normalised on $W_{mb}$.
\end{proof}

Therefore
we may define
\[
\mu(g,p)=g\cdot p=a_{y}\cdot m_{x}=(am)_{y}
\]
and
\[
\eta(p,h)=p\cdot h=m_{x}\cdot b_{z}=(mb)_{x}\;.
\]
We have to show that this two definitions are independent of the choice of
$a$, $m$ and $b$.
To this end, observe that the left $A$-action on $\cM$,
as a map
$A \ten_{A_{0}}\cM \ra \cM$,
is a homomorphism of left $A_{0}$-coalgebras,
and induces
a homomorphism of left $(A_{0})_{y}$-coalgebras
\[
A_{y}\ten_{A_{0}}\cM\cong (A \ten_{A_{0}}\cM)_{y} \ra \cM_{y}\;.
\]
By Lemma \ref{lemmagrouplikeelements}
we know that $[g,p]\in G(A_{y}\ten_{A_{0}}\cM)$,
so its image in $\cM_{y}$ is a grouplike element.
A direct computation shows that this image is exactly $\mu(g,p)$,
which shows that $\mu$ is well defined as a map
\[
\mu\!:G\times^{s,\pi}_{G_{0}}P\ra P\;.
\]
Similarly, the right $B$-action
$\cM \ten_{B_{0}} B\ra \cM$ on $\cM$
is a homomorphism of left
$A_{0}$-coalgebras, and gives a homomorphism of left
$(A_{0})_{x}$-coalgebras
\[
\cM_{x}\ten_{B_{0}}B \cong (\cM \ten_{B_{0}} B)_{x} \ra \cM_{x}\;.
\]
Again by Lemma \ref{lemmagrouplikeelements}
we know that $[p,h]\in G(\cM_{x}\ten_{B_{0}}B)$,
so its image in $\cM_{x}$ is a grouplike element,
equal to $\eta(p,h)$. This shows that
$\eta$ is well defined as a map
\[
\eta\!:P\times^{\phi,t}_{H_{0}}H \ra P\;.
\]

\begin{prop}
The map $\mu$ is a left action of $G$ on $P$ along
$\pi$.
\end{prop}

\begin{proof}
Suppose that $m\in G_{w}(\cM)$ is normalised on an open
neighbourhood $W_{m}\subset G_{0}$ of a point $x\in M$.
Let $a,a'\in G_{w}^{S}(A)$ be normalised on open subsets $W_{a}$
respectively $W_{a'}$ of $G_{0}$, and let $y\in W_{a}$ and
$y'\in W_{a'}$ be such that $a_{y}$ is an arrow from $x$ to $y$ and
$a'_{y'}$ is an arrow from $y$ to $y'$.

Straight from the definition of $\mu$ it follows
$\pi(a_{y}\cdot m_{x})=y$, thus $\mu$ acts along the
map $\pi$. Furthermore, we have
\[
(a'_{y'}a_{y})\cdot m_{x}=(a'a)_{y'}\cdot m_{x}=(a'am)_{y'}
=a'_{y'}\cdot(am)_{y}=a'_{y'}\cdot(a_{y}\cdot m_{x})\;.
\]
The unit arrow $1_{x}\in G$ can be represented by a
smooth function $u_{0}\in\Cc(G_{0})$ satisfying $(u_{0})_{x}=1$, thus
\[
1_{x}\cdot m_{x}=(u_{0})_{x}\cdot
m_{x}=(u_{0}m)_{x}=(u_{0})_{x}m_{x}=m_{x}\;.
\]

To prove that $\mu$ is smooth,
first observe that the map $w\!:G\times_{G_{0}}P\ra G_{0}$,
given by $w=s\com \pr_{1}=\pi\com\pr_{2}$,
is a local diffeomorphism. The
neighbourhood $W_{1}=a_{W_{a}}\times_{G_{0}}m_{W_{m}}$ of the
point $(a_{y},m_{x})\in G\times_{G_{0}}P$
is mapped by $w$
diffeomorphically onto $V_{W_{a},a}\cap W_{m}$. Furthermore, the
neighbourhood $W_{2}=(am)_{W_{am}}$ of the point
$\mu(a_{y},m_{x})\in P$ is mapped diffeomorphically onto
$W_{am}=\tau_{W_{a},a}(V_{W_{a},a}\cap W_{m})$ by
$\pi$. We can locally express the map $\mu$ as
$\mu|_{W_{1}}=(\pi|_{W_{2}})^{-1}\com\tau_{W_{a},a}\com
w|_{W_{1}}$, which shows that the action is smooth.
\end{proof}

\begin{prop}
The map $\eta$ is a right action of $H$ on $P$ along
$\phi$.
\end{prop}

\begin{proof}
Let $b,b'\in G_{w}^{S}(B)$ be normalised on open subsets
$W_{b}$ respectively $W_{b'}$ of $H_{0}$,
and suppose that $z\in W_{b}$ and $z'\in W_{b'}$ are such that
$b_{z}$ is an arrow from $z'$ to $z$
and $b'_{z'}$ is an arrow from $z''$ to $z'$.
Furthermore, assume that $m\in G_{w}(\cM)$ is normalised on an
open neighbourhood $W_{m}$ of a point $x$ in $G_{0}$ such that
$\phi(m_{x})=z$. Choose
$m'\in\cM$ with $\cm(m)=m\ten m'$.
Since
$\cu(m'v'_{0})|_{W}=\cu(mv'_{0})|_{W}$ and $mv'_{0}=\cu(m'v'_{0})m$
for any $v'_{0}\in\Cc(H_{0})$, it follows that
\begin{align*}
v_{0}(\phi(m_{x}\cdot b_{z}))&= \cu(mbv_{0})(x)\\
                             &= \cu(m(v_{0}\com\tau_{W_{b},b}^{-1})b)(x)\\
                             &= \cu(\cu(m'(v_{0}\com\tau_{W_{b},b}^{-1}))mb)(x)\\
                             &= \cu(m'(v_{0}\com\tau_{W_{b},b}^{-1}))(x)\cu(mb)(x)\\
                             &= \cu(m'(v_{0}\com\tau_{W_{b},b}^{-1}))(x)\\
                             &= \cu(m(v_{0}\com\tau_{W_{b},b}^{-1}))(x)\\
                             &= (v_{0}\com\tau_{W_{b},b}^{-1})(\phi(m_{x}))\\
                             &= v_{0}(z')
\end{align*}
for arbitrary $v_{0}\in\Cc(V_{W_{b},b})$. If $\phi(m_{x}\cdot b_{z})$
and $z'\in V_{W_{b},b}$ were different points of $H_{0}$, we could choose
$v_{0}\in\Cc(V_{W_{b},b})$ such that
$v_{0}(\phi(m_{x}\cdot b_{z}))\neq v_{0}(z')$.
The above calculation thus shows that
$\phi(m_{x}\cdot b_{z})=z'$.

Next we  have
\[
(m_{x}\cdot b_{z})\cdot
b'_{z'}=(mb)_{x}\cdot
b'_{z'}=(mbb')_{x}=m_{x}\cdot(bb')_{z}=m_{x}\cdot(b_{z}b'_{z'})\;.
\]
If we represent the identity arrow $1_{z}$ by $v_{0}\in\Cc(H_{0})$
with $(v_{0})_{z}=1$, we get
\[
m_{x}\cdot
1_{z}=m_{x}(v_{0})_{z}=m_{x}v_{0}=(v_{0}\com\phi_{W_{m},m})_{x}m_{x}=m_{x}.
\]

Finally, we show that $\eta$
is smooth.
Note that the projection
$\pr_{1}\!:P\times_{H_{0}}H\ra P$ is a local
diffeomorphism.
The neighbourhood
$W_{1}=m_{W_{m}}\times_{H_{0}}b_{W_{b}}$ of the point
$(m_{x},b_{z})\in P\times_{H_{0}}H$
is mapped by $\pi\com \pr_{1}$
diffeomorphically onto
$W_{mb}=\pi(m_{W_{m}}\cap\phi^{-1}(W_{b}))$. Similarly, the
neighbourhood $W_{2}=(mb)_{W_{mb}}$ of the point
$\eta(m_{x},b_{z})\in P$ is mapped by $\pi|_{W_{2}}$
diffeomorphically onto $W_{mb}$. We can locally express $\eta$ as
$\eta|_{W_{1}}=(\pi|_{W_{2}})^{-1}\com(\pi\com \pr_{1})|_{W_{1}}$
and conclude that $\eta$ is smooth.
\end{proof}

\subsection{Principalness of $P$}

Finally, we need to show that $P$ constructed above
is indeed a principal $H$-bundle over $G$.

\begin{prop}
The manifold $P$, with the actions $\mu$ and $\eta$, is
a principal $H$-bundle over $G$.
\end{prop}

\begin{proof}
Let $g\in G(x,y)$, $p\in P$ and $h\in H(z',z)$
be such that $\pi(p)=x$ and $\phi(p)=z$.
Choose $a\in G_w^S(A)$, $m\in G_w(\cM)$,
and $b\in G_{w}^{S}(B)$ such that
$a$ is normalised on an open neighbourhood $W_{a}\subset G_{0}$ of $y$,
$m$ is normalised on an open neighbourhood $W_{m}\subset G_{0}$ of $x$
and
$b$ is normalised on an open neighbourhood $W_{b}\subset H_{0}$ of $z$
with $g=a_{y}$, $p=m_{x}$ and $h=b_{z}$.

(i) Straight from the definition of the action of $H$ on $P$ it
follows that $\pi(m_{x}\cdot b_{z})=\pi(m_{x})$, which proves that
$H$ acts along the fibers of the map $\pi$. Next, choose arbitrary
$v_{0}\in\Cc(H_{0})$ and $u_{0}\in\Cc(W_{a})$ with
$(u_{0})_{y}=1$. Then
\begin{align*}
v_{0}(\phi(a_{y}\cdot m_{x}))
  & = \cu(amv_{0})(y)=u_{0}(y)\cu(amv_{0})(y)\\
  & = \cu(u_{0}a\cu(mv_{0}))(y)=\cu(a(u_{0}\com\tau_{W_{a},a})\cu(mv_{0}))(y)\\
  & = \cu(u_{0}(\cu(mv_{0})\com\tau_{W_{a},a}^{-1})a)(y)\\
  & = u_{0}(y)(\cu(mv_{0})\com\tau_{W_{a},a}^{-1})(y)\cu(a)(y)\\
  & = \cu(mv_{0})(x)\cu(a)(y)=\cu(mv_{0})(x)\\
  & = v_{0}(\phi(m_{x})).
\end{align*}
The function $u_{0}\in\Cc(W_{a})$ was used to ensure that
$(u_{0}\com\tau_{W_{a},a})\cu(mv_{0})\in\Cc(V_{W_{a},a})$. Since
$v_{0}\in\Cc(H_{0})$ was arbitrary, we conclude that
$\phi(a_{y}\cdot m_{x})=\phi(m_{x})$, which shows that $G$ acts
along the fibers of the map $\phi$.

(ii) Both actions commute as a result of
\[
(a_{y}\cdot m_{x})\cdot b_{z}=(am)_{y}\cdot
b_{z}=(amb)_{y}=a_{y}\cdot(mb)_{x}=a_{y}\cdot(m_{x}\cdot b_{z})\;.
\]

(iii) The map $\pi$ is surjective because $\cu$ is surjective.
Finally, we have to show that $(\pr_{1},\eta)\!:P\times_{H_{0}}H\ra
P\times_{G_{0}}P$ is a diffeomorphism.
To see that it is a bijection, it is sufficient to
show that it restricts to a bijection
between the corresponding
fibers over any $x\in G_{0}$, that is from
$\bigcup_{p\in\pi^{-1}(x)}(\{p\}\times t^{-1}(\phi(p)))
\subset P\times_{H_{0}}H$
to
$\pi^{-1}(x)\times\pi^{-1}(x)\subset P\times_{G_{0}}P$.
By Lemma \ref{lemmadeltabar} we know that
$\overline{\cm}_{x}\!:\cM_{x}\ten_{B_{0}}B\ra \cM_{x}\ten_{A_{0}}
\cM_{x}$ is an isomorphism of coalgebras over $(A_{0})_{x}$, so it
restricts to a bijection between
$G(\cM_{x}\ten_{B_{0}}B)$ and $G(\cM_{x}\ten_{A_{0}}\cM_{x})$.
By Lemma \ref{lemmagrouplikeelements},
applied to $\cN=B$, we may
identify $G(\cM_{x}\ten_{B_{0}}B)$ with
$\bigcup_{p\in\pi^{-1}(x)}(\{p\}\times t^{-1}(\phi(p)))
\subset P\times_{H_{0}}H$
and
$G(\cM_{x}\ten_{A_{0}}\cM_{x})$ with
$G(\cM_{x})\times G(\cM_{x})=\pi^{-1}(x)\times\pi^{-1}(x)$.
Observe that the bijection
\[ \bigcup_{p\in\pi^{-1}(x)}(\{p\}\times t^{-1}(\phi(p)))\ra
   \pi^{-1}(x)\times\pi^{-1}(x) \]
given by $\overline{\cm}_{x}$ is in fact of the form
$(p,h)\mapsto(p,p\cdot h)$, thus equal to the restriction of
$(\pr_{1},\eta)\!:P\times_{H_{0}}H\ra P\times_{G_{0}}P$ to the fiber
over $x$. This proves that $(\pr_{1},\eta)\!:P\times_{H_{0}}H\ra
P\times_{G_{0}}P$ is a bijection. Furthermore, the map
$(\pr_{1},\eta)$ is a local diffeomorphism because $H$ is
\'{e}tale, and therefore it is a diffeomorphism.
\end{proof}

\section{Equivalence of the Morita categories}

In this section we state and prove the main
result of this paper:

\begin{theo}\label{maintheo}
The functor
$\Cc\!:\EtGPD\ra\LgHoALGD$
is an equivalence between
the Morita category of \'{e}tale Lie groupoids
and the Morita category
of locally grouplike Hopf algebroids.
\end{theo}

Before we give the proof, let us start with two lemmas.
Let
$(A,G_{0})$ and $(B,H_{0})$ be locally grouplike Hopf algebroids
with the associated spectral \'{e}tale groupoids
$G$ respectively $H$.
We denote by $\cM_{\ast}$ the principal $H$-bundle $P$ over $G$
associated to a locally grouplike principal $A$-$B$-bimodule
$\cM$, as in Section \ref{section3}.
Suppose that
$\theta\!:\cM\ra\cM'$ is a homomorphism of principal
$A$-$B$-bimodules. In particular, $\theta\!:\cM\ra\cM'$ is a
homomorphism of coalgebras over $A_{0}=\Cc(G_{0})$,
so we have the
associated map
$\Esp(\theta)\!:
 \Esp(\cM)\ra\Esp(\cM')$ of spectral sheaves over $G_{0}$
\cite{Mrc02,Mrc0x}.
We can describe the map $\Esp(\theta)$ as follows:
If a point $p\in G(\cM_{x})=\Esp(\cM)_{x}$ is represented
by an element $m\in G_{w}(\cM)$ normalised at $x\in G_{0}$
(i.e.\ $p=m_{x}$), then
$\theta(m)\in G_{w}(\cM')$ is normalised at $x$ as well and
\[
\Esp(\theta)(p)=(\theta(m))_{x}\in G(\cM'_{x})=\Esp(\cM')_{x}\;.
\]
Since $\cM_{\ast}=\Esp(\cM)$ and $\cM'_{\ast}=\Esp(\cM')$
as sheaves over $M$, we shall write
$\theta_{\ast}=\Esp(\theta)\!:\cM_{\ast}\ra\cM'_{\ast}$.

\begin{lem}\label{lemfunctor}
The map $\theta_{\ast}\!:\cM_{\ast}\ra\cM'_{\ast}$ is an equivariant map of
principal $H$-bundles over $G$.
\end{lem}

\begin{proof}
Suppose that $m\in G_{w}(\cM)$ is normalised at $x\in G_{0}$,
$a\in G_{w}^{S}(A)$ is normalised at $y\in G_{0}$
and $b\in G_{w}^{S}(B)$ is normalised at $\phi(m_{x})=z$,
such that $a_{y}\in G(x,y)$ and $b_{z}\in H(z',z)$.
Straight from the definition of $\theta_{\ast}$ it follows
that $\pi'\com\theta_{\ast}=\pi$. Furthermore, for any
$v_{0}\in\Cc(H_{0})$ we have
\[
v_{0}(\phi'(\theta_{\ast}(m_{x})))=\cu'(\theta(m)v_{0})(x)
=\cu'(\theta(mv_{0}))(x)
=\cu(mv_{0})(x)=v_{0}(\phi(m_{x}))\;,
\]
which shows that $\phi'\com\theta_{\ast}=\phi$.
Next, the equalities
\[
\theta_{\ast}(a_{y}\cdot m_{x})=\theta_{\ast}((am)_{y})=\theta(am)_{y}=
(a\theta(m))_{y}=a_{y}\cdot\theta(m)_{x}=a_{y}\cdot\theta_{\ast}(m_{x})
\]
and
\[
\theta_{\ast}(m_{x}\cdot b_{z})=\theta_{\ast}((mb)_{x})=\theta(mb)_{x}=
(\theta(m)b)_{x}=\theta(m)_{x}\cdot b_{z}=\theta_{\ast}(m_{x})\cdot
b_{z}
\]
show that $\theta_{\ast}$ is equivariant.
\end{proof}

For any principal $H$-bundle $P$ over $G$ there is a natural
isomorphism $\Phi\!:P\ra\Esp(\Cc(P))=\Cc(P)_{\ast}$
of sheaves over $G_{0}$
\cite{Mrc02,Mrc0x}, given by
\[
\Phi(p)=f_{\pi(p)}\;,
\]
where $f\in\ G_{w}(\Cc(P))$ with $f_{p}=1$.
With respect to the natural isomorphisms
$G\cong\Gsp(\Cc(G))$ and $H\cong\Gsp(\Cc(H))$, given by the same
formula as $\Phi$, we may regard $\Cc(P)_{\ast}$ as a principal
$H$-bundle over $G$.

\begin{lem}\label{lembsp}
The map $\Phi\!:P\ra\Cc(P)_{\ast}$ is an isomorphism of principal
$H$-bundles over $G$.
\end{lem}
\begin{proof}
Let $g\in G(x,y)$, $h\in H(z',z)$ and
$p\in P$ with $\pi(p)=x$ and $\phi(p)=z$. Suppose that
$u\in G_{w}^{S}(\Cc(G))$, $v\in G_{w}^{S}(\Cc(H))$ and $f\in
G_{w}(\Cc(P))$ satisfy $u_{g}=1$,
$v_{h}=1$ and $f_{p}=1$.
For any
$v_{0}\in\Cc(H_{0})$ we have
$(fv_{0})(p)=f(p)v_{0}(\phi(p))=v_{0}(\phi(p))$,
by the equation (\ref{eqrightaction}).
Since $fv_{0}\in G_{w}(\Cc(P))$ we have
$(fv_{0})(p)=\cu(fv_{0})(\pi(p))$
and hence
\[
v_{0}(\phi(p))=\cu(fv_{0})(\pi(p))=v_{0}(\phi(f_{\pi(p)}))
=v_{0}(\phi(\Phi(p)))\;,
\]
which shows that $\phi\com\Phi=\phi$.

By Lemma \ref{lem1} we have $uf\in G_{w}(\Cc(P))$, $fv\in
G_{w}(\Cc(P))$, $\cu(uf)_{y}=1$ and $\cu(fv)_{x}=1$.
Furthermore, we have the equalities $(uf)_{g\cdot p}=1$
respectively $(fv)_{p\cdot h}=1$ (by the equations
(\ref{eqleftaction}) and (\ref{eqrightaction})), so
\[
\Phi(g\cdot p)=(uf)_{y}=u_{y}\cdot f_{x}
=g\cdot\Phi(p)
\]
and
\[
\Phi(p\cdot h)=(fv)_{x}=f_{x}\cdot v_{z}
=\Phi(p)\cdot h\;,
\]
which shows that $\Phi$ is equivariant.
\end{proof}

\begin{proof}[Proof of Theorem \ref{maintheo}]
(i)
Any locally grouplike Hopf algebroid $(A,M)$
is isomorphic to the locally grouplike Hopf algebroid
$(\Cc(\Gsp(A)),M)$ \cite{Mrc07}, so $\Cc$ is
essentially surjective.

(ii) Let $G$ and $H$ be \'{e}tale Lie groupoids and let $\cM$ be a
locally grouplike principal $\Cc(G)$-$\Cc(H)$-bimodule. The map $
\Psi\!:\Cc(\cM_{\ast})\ra\cM$, given by
\[
\Psi(\sum_{i=1}^{k}f_{i}\com
\pi|_{(m_{i})_{W_{m_{i}}}})=\sum_{i=1}^{k}f_{i}m_{i}\;,
\]
is an isomorphism of coalgebras over $\Cc(G_{0})$
(\cite[Theorem 2.5]{Mrc02}, see also \cite{Mrc0x}),
where $m_{i}\in G_{w}(\cM)$,
$W_{m_{i}}$ is an open subset of $G_{0}$ such that
$m_{i}$ is normalised on $W_{m_{i}}$
and $f_{i}\in\Cc(W_{m_{i}})$, for any $i=1,\ldots,k$.

We will show that
$\Psi(uf)=u\Psi(f)$ and $\Psi(fv)=\Psi(f)v$ for any
$u\in\Cc(G)$, $f\in\Cc(\cM_{\ast})$ and $v\in\Cc(H)$.
We can assume that
$f=f_{0}\com\pi|_{U}$, where $U=m_{W_{m}}$
for some $m\in G_{w}(\cM)$ normalized on an
open subset $W_{m}\subset G_{0}$
and $f_{0}\in\Cc(W_{m})$.
Furthermore, we may assume that
$u=u_{0}\com t|_{U'_{u}}$ and $v=v_{0}\com t|_{U'_{v}}$, where
$U'_{u}\subset G$ respectively $U'_{v}\subset H$ are bisections,
$u_{0}\in\Cc(t(U'_{u}))$ and $v_{0}\in\Cc(t(U'_{v}))$.
Choose elements
$a\in \Cc(U'_{u})\subset G_{w}^{S}(\Cc(G))$, normalised on an open
neighbourhood $W_{a}\subset t(U'_{u})$ of $\supp(u_{0})$, and
$b\in \Cc(U'_{v})\subset G_{w}^{S}(\Cc(H))$,
normalised on an open neighbourhood
$W_{b}\subset t(U'_{v})$ of $\supp(v_{0})$.
Write
$U_{u}=t^{-1}(W_{a})\cap U'_{u}$
and
$U_{v}=t^{-1}(W_{b})\cap U'_{v}$.
By the equations (\ref{eqleftaction}) and
(\ref{eqrightaction}) we get
\[
uf=(u_{0}(f_{0}\com\tau_{W_{a},a}^{-1}))\com\pi|_{\mu(U_{u}\times_{G_{0}}U)}
\]
and
\[
fv=(f_{0}(v_{0}\com\phi_{W_{m},m}))\com\pi|_{\eta(U\times_{H_{0}}U_{v})}\;.
\]
Denote $W_{am}=\tau_{W_{a},a}(V_{W_{a},a}\cap W_{m})$ and
$W_{mb}=\pi(m_{W_{m}}\cap\phi^{-1}(W_{b}))$. The elements
$am,mb\in G_{w}(\cM)$ are, by Lemma \ref{lem1}, normalised on the
supports of the functions
$u_{0}(f_{0}\com\tau_{W_{a},a}^{-1})\in\Cc(W_{am})$ respectively
$f_{0}(v_{0}\com\phi_{W_{m},m})\in\Cc(W_{mb})$. Moreover,
$\mu(U_{u}\times_{G_{0}}U)=(am)_{W_{am}}$ and
$\eta(U\times_{H_{0}}U_{v})=(mb)_{W_{mb}}$, which implies
\[
\Psi(uf)=u_{0}(f_{0}\com\tau_{W_{a},a}^{-1})am
\]
and
\[
\Psi(fv)=f_{0}(v_{0}\com\phi_{W_{m},m})mb\;.
\]
On the other hand we have
$u\Psi(f)=uf_{0}m=u_{0}af_{0}m=u_{0}(f_{0}\com\tau_{W_{a},a}^{-1})am$
respectively
$\Psi(f)v=f_{0}mv=f_{0}mv_{0}b=f_{0}(v_{0}\com\phi_{W_{m},m})mb$,
which shows that $\Psi$ is an isomorphism of locally
grouplike principal bimodules. This proves that the functor $\Cc$
is full.

(iii) Finally, we show that $\Cc$ is faithful. Choose principal
$H$-bundles $P$ and $P'$ over $G$, and suppose that there exists an
isomorphism $\theta\!:\Cc(P)\ra\Cc(P')$ of locally grouplike
principal $\Cc(G)$-$\Cc(H)$-bimodules. The principal bundles
$\Cc(P)_{\ast}$
and $\Cc(P')_{\ast}$ are isomorphic by Lemma \ref{lemfunctor},
so $P$ and $P'$ are isomorphic as well by Lemma \ref{lembsp} and
therefore represent the same morphism from $G$ to $H$ in the
category $\EtGPD$.
\end{proof}

\end{document}